\documentclass[11pt]{article}
\usepackage{amsmath, amscd, amssymb, latexsym,
    epsfig, color, amsthm}
\numberwithin{equation}{section}
\def\proof{\smallskip\noindent {\it Proof: \ }}
\def\endproof{\hfill$\square$\medskip}
\newtheorem{theorem}{Theorem}[section]

\newtheorem{corollary}[theorem]{Corollary}

\newtheorem{lemma}[theorem]{Lemma}

\theoremstyle{definition}

\newtheorem{remark}[theorem]{Remark}

\def\R{\mathbb{R}}
\def\Z{\mathbb{Z}}
\def\S{\mathbb{S}}

\def\P{\mathcal{P}}
\def\F{\mathcal{F}}

\def\B{\mathcal{B}}
\def\S{\mathbb{S}}

\DeclareMathOperator\conv{\mathrm{conv}}

\begin{document}
\title{\vspace{-1ex}
Explicit constructions of centrally symmetric
$k$-neighborly polytopes and large strictly antipodal
sets}

\author{Alexander Barvinok
\thanks{Research of the first and second authors
is partially supported by NSF grant DMS-0856640}\\
\small  Department of Mathematics \\[-0.8ex]
\small University of Michigan\\[-0.8ex]
\small Ann Arbor, MI 48109-1043, USA\\[-0.8ex]
\small \texttt{barvinok@umich.edu}
\and Seung Jin Lee \\
\small  Department of Mathematics \\[-0.8ex]
\small University of Michigan\\[-0.8ex]
\small Ann Arbor, MI 48109-1043, USA\\[-0.8ex]
\small \texttt{lsjin@umich.edu}
\and Isabella Novik
\thanks{Research of the third author is partially
supported by NSF grant DMS-1069298}\\
\small Department of Mathematics, Box 354350\\[-0.8ex]
\small University of Washington\\[-0.8ex]
\small Seattle, WA 98195-4350, USA\\[-0.8ex]
\small \texttt{novik@math.washington.edu}
}

\maketitle
\begin{abstract}
We present explicit constructions of
centrally symmetric $2$-neighborly $d$-dimensional
polytopes with about $3^{d/2}\approx (1.73)^d$ vertices
and of centrally symmetric $k$-neighborly $d$-polytopes
with about $2^{{3d}/{20k^2 2^k}}$ vertices. Using this
result, we construct for a fixed $k\geq 2$ and arbitrarily
large $d$ and $N$, a centrally symmetric $d$-polytope
with $N$ vertices that has at least
$\left(1-k^2\cdot(\gamma_k)^d\right)\binom{N}{k}$
faces of dimension $k-1$, where
$\gamma_2=1/\sqrt{3}\approx 0.58$ and
$\gamma_k = 2^{-3/{20k^2 2^k}}$ for $k\geq 3$.
Another application is a construction of a set
of $3^{\lfloor d/2 -1\rfloor}-1$ points in $\R^d$ every two
of which are strictly antipodal as well as a construction
of an $n$-point set (for an arbitrarily large $n$)
in $\R^d$ with many pairs of strictly
antipodal points. The two latter results significantly
improve the previous bounds by Talata, and Makai and Martini,
respectively.
\end{abstract}

\section{Introduction}

\subsection{Cs neighborliness}
What is the maximum number of $k$-dimensional faces that a
centrally symmetric $d$-dimensional polytope with $N$ vertices
can have? While the answer in the class of all polytopes
is classic by now \cite{McM-UBC}, very little is known
in the centrally symmetric case. Here we present several
constructions that significantly improve existing lower
bounds on this number.

Recall that a {\em polytope} is the convex hull of a
set of finitely many points in $\R^d$. The dimension of
a polytope $P$ is the dimension of its affine hull.
A polytope $P \subset \R^d$ is
{\it centrally symmetric} (cs, for short) if $P =-P$.
A cs polytope $P$ is {\em $k$-neighborly} if
every set of $k$ vertices of $P$ no two of which are
antipodes forms the vertex set of a $(k-1)$-face of $P$.

It was proved in \cite{LinNov} that a cs $2$-neighborly
$d$-dimensional polytope cannot have more than $2^d$ vertices.
On the other hand, a construction from \cite{BLN-many}
showed that there exist such polytopes with about
$3^{d/4}\approx(1.316)^d$ vertices.
In Theorem \ref{2-neighb}(1) we present a construction
of a cs $2$-neighborly $d$-polytope with about
$3^{d/2}\approx(1.73)^d$ vertices.

More generally, it was verified in \cite{BN} that a cs
$d$-dimensional polytope with $N$ vertices cannot have more
than $\left(1-0.5^{d}\right) \frac{N^2}{2}$ edges. However,
a construction from \cite{BLN-many} produced cs
$d$-dimensional polytopes with $N$ vertices and about
$\left(1-3^{-d/4}\right) \frac{N^2}{2}
\approx \left(1-0.77^{d}\right) \frac{N^2}{2}$ edges.
In Theorem \ref{2-neighb}(2), we improve this bound by
constructing a cs
$d$-dimensional polytope with $N$ vertices (for an
arbitrarily large $N$) and at least
$\left(1-3^{-\lfloor d/2-1\rfloor}\right) \binom{N}{2}
\approx \left(1-0.58^{d}\right) \frac{N^2}{2}$ edges.

For higher-dimensional faces even less is known.
It follows from the results of \cite{LinNov}
that no cs $k$-neighborly $d$-polytope can have more than
$\lfloor d\cdot 2^{{C}d/k} \rfloor$ vertices,
where ${C}>0$ is some absolute constant.
At the same time, papers \cite{LinNov, RV} used a randomized
construction to prove existence of
$k$-neighborly cs $d$-dimensional polytopes
with $\lfloor d\cdot 2^{cd/k} \rfloor$ vertices for
some absolute constant $c>0$.
However, for  $k>2$ no deterministic construction of a $d$-dimensional 
$k$-neighborly cs polytope with $2^{\Omega(d)}$ vertices is known.
In Theorem \ref{k-neighb} and Remark \ref{numerics}
we present a deterministic construction of a cs $k$-neighborly
$d$-polytope with at least $2^{c_k d}$ vertices where $c_k=3/20k^2 2^k$.
We then use this result in Corollary~\ref{many-k}
to construct for a fixed $k$ and
arbitrarily large $N$ and $d$, a cs $d$-polytope with $N$
vertices that has a record number of $(k-1)$-dimensional faces.
Our construction relies on the notion of $k$-independent
families \cite{Alon, FLL} (see also \cite{AMS}).

Through Gale duality $m$-dimensional subspaces of $\R^N$ 
correspond to $(N-m)$-dimensional cs polytopes with $2N$ vertices. 
If the subspace is ``almost Euclidean'' (meaning that the ratio 
of the $\ell^1$ and $\ell^2$ norms of
nonzero vectors of the subspace remains within certain bounds, 
see \cite{LinNov} for technical details), then the corresponding polytope 
turns out to be $k$-neighborly. 
Despite considerable efforts, see for example \cite{Indyk}, no
explicit constructions of ``almost Euclidean'' subspaces is known 
for $m$ anywhere close to $N$. Our polytopes give rise 
to subspaces of $\R^N$ of codimension $O(\log N)$ and it would be interesting 
to find out if the resulting subspaces are indeed ``almost Euclidean''.

\subsection{Antipodal points}
Our results on cs polytopes provide new bounds on
several problems related to strict antipodality.
Let $X\subset\R^d$ be a set that affinely spans $\R^d$.
A pair of points $u,v\in X$
is called {\em strictly antipodal} if there exist two distinct
parallel hyperplanes $H$ and $H'$ such that
$X\cap H=\{u\}$, $X\cap H'=\{v\}$, and
$X$ lies in the slab between $H$ and $H'$.
Denote by $A'(d)$ the maximum size of a set $X\subset\R^d$
having the property that every pair of points of $X$
is strictly antipodal, by $A'_d(Y)$ the number of
strictly antipodal pairs of a given set $Y$,
and by $A'_d(n)$ the maximum size of $A'_d(Y)$
taken over all $n$-element subsets $Y$ of $\R^d$.
(Our notation follows the recent survey paper \cite{MartSolt}.)

The notion of strict antipodality was introduced in 1962
by Danzer and Gr\"unbaum \cite{DanzGr} who verified  that
$2d-1\leq A'(d)\leq 2^d$ and conjectured that $A'(d)=2d-1$.
However, twenty years later, Erd\H{o}s and
F\"uredi \cite{ErdFur} used a probabilistic argument to prove
that $A'(d)$ is exponential in $d$.  
 Their result was improved by
Talata (see \cite[Lemma 9.11.2]{Bor}) who
found an explicit construction showing that for $d\geq 3$,
\[
A'(d)\geq \lfloor (\sqrt[3]{3})^d/3\rfloor.
\]
Talata also announced that $(\sqrt[3]{3})^d/3$ in the above
formula can be replaced with $(\sqrt[4]{5})^d/4$. 
(It is worth remarking that Erd\H{o}s and F\"uredi  established 
existence of an {\em acute set} in $\R^d$ that has 
an exponential size in $d$. As every acute set has the property that
all of its pairs of vertices are strictly antipodal, their result
implied an exponential lower bound on $A'(d)$. A significant improvement
of the  Erd\H{o}s--F\"uredi bound on the maximum size of an acute set in 
$\R^d$ was recently found by Harangi \cite{Ha}.)

Regarding the value of $A'_d(n)$,
Makai and Martini \cite{MakMart} showed that for $d\geq 4$,
\[
\left(1-\frac{\mbox{const}}{(1.0044)^d}\right)\frac{n^2}{2}-O(1)
\leq A'_d(n)\leq \left(1-\frac{1}{2^{d}-1}\right)\frac{n^2}{2}.
\]

Here we observe that an appropriately chosen half of the
vertex set of a cs $d$-polytope with many edges
has a large number of
strictly antipodal pairs of points. Consequently,
our construction of cs $d$-polytopes with many edges
implies --- see Theorem \ref{A'_d} --- that
$$A'(d)\geq 3^{\lfloor d/2-1\rfloor}-1 \quad \mbox{and}
\quad
A'_d(n)\geq
\left(1-\frac{1}{3^{\lfloor d/2-1\rfloor}-1}\right)
\frac{n^2}{2} -O(n)
\quad \mbox{for all } d\geq 4.
$$

\bigskip
The rest of the paper is structured as follows.
In Section 2 we review several facts and definitions
related to the symmetric moment curve. In Section 3,
we present our construction of a cs 2-neighborly
$d$-polytope with many vertices as well as that of a cs
$d$-polytope with arbitrarily
many vertices and a record number of edges. Section 4
is devoted to applications of these results to problems
on strict antipodality. Finally, in Section 5 we provide
a deterministic construction of a cs $k$-neighborly
$d$-polytope and of a cs $d$-polytope with arbitrarily
many vertices and a record number of $(k-1)$-faces.

\section{The symmetric moment curve}
In this section we collect several definitions and results
needed for the proofs. We start with the notion of
the symmetric moment curve on which all our
constructions are based. The {\em symmetric moment curve}
$U_k : \R\longrightarrow \R^{2k}$ is defined by
\begin{equation}  \label{U_k}
U_k(t)=\left(\cos t,\ \sin t,\ \cos 3t,\ \sin3t,\ \ldots,\
\cos (2k-1)t,\ \sin(2k-1)t \right).
\end{equation}
 Since
$$U_k(t) = U_k(t + 2\pi) \quad \text{for all} \quad t,$$
from this point on, we consider $U_k(t)$ to be defined on the
unit circle $\S =\R /2 \pi \Z.$
We note that $t$ and $t+\pi$ form a pair of antipodes
for all $t \in \S$ and that
$$U_k(t+\pi)=-U_k(t) \quad \text{for all} \quad t \in \S.$$

The value of an affine function 
$A: \R^{2k}\longrightarrow\R$ on the symmetric
moment curve $U_k$ is represented by a trigonometric polynomial
of degree at most $2k-1$ that has the following form
\[f(t)=c+\sum_{j=1}^k a_j\cos(2j-1)t+\sum_{j=1}^k b_j\sin(2j-1)t,
\quad \mbox{where } a_j,b_j,c\in \R. 
\]
Starting with any trigonometric polynomial $f:\S\longrightarrow \R$,
$f(t)=c+\sum_{j=1}^d a_j\cos(jt)+\sum_{j=1}^d b_j\sin(jt)$
of degree at most $d$
and substituting $z=e^{it}$ gives rise to a complex polynomial
\begin{equation}  \label{P(f)}
\P(f)(z):=z^{d}\left(c+\sum_{j=1}^d a_j \frac{z^{j}+z^{-j}}{2}
+\sum_{j=1}^d b_j \frac{z^{j}-z^{-j}}{2i} \right).
\end{equation}
This polynomial has degree at most $2d$,
it is self-inversive (that is, the coefficient of $z^j$
is conjugate to that of $z^{2d-j}$), and $t^\ast\in\S$ is a root of
$f(t)$ if and only if $e^{it^\ast}$ is a root of $\P(f)(z)$
(see \cite{BN} and \cite{BLN-neighb} for more details).  In particular,  $f(t)$ cannot have
more than $2d$ roots (counted with multiplicities).

The following result concerning the convex hull of
the symmetric moment curve was proved in \cite{BLN-neighb}. 
In what follows we talk about exposed faces, that is, intersections
of convex bodies with supporting affine hyperplanes.
\begin{theorem}   \label{pi/2}
Let $\B_k \subset {\R}^{2k}$,
$$\B_k=\conv\left(U_k(t): \quad t \in \S\right),$$
be the convex hull of the symmetric moment curve. Then
for every positive integer $k$ there exists a number
$$\frac{\pi}{2} \ < \ \alpha_k \ < \ \pi$$
such that for an arbitrary open arc $\Gamma \subset {\S}$
of length $\alpha_k$ and arbitrary distinct $n \leq k$ points
$t_1, \ldots, t_n \in \Gamma$, the set
$$\conv\left(U_k\left(t_1\right), \ldots,
U_k\left(t_n\right)\right)$$
is a face of $\B_k$.
\end{theorem}

\noindent For $k=2$ with $\alpha_2=2\pi/3$ this result
is due to Smilansky \cite{Sm85}.

We also frequently use the following well-known fact
about polytopes: if
$T: {\R}^{d'} \longrightarrow {\R}^{d''}$ is a linear
transformation and $P \subset {\R}^{d'}$ is a polytope,
then $Q=T(P)$ is also a polytope and for every face $F$
of $Q$ the inverse image of $F$,
$$T^{-1}(F)=\left\{x \in P: \quad T(x) \in F \right\},$$
is a face of $P$; this face is the convex hull of the vertices
of $P$ mapped by $T$ into vertices of $F$.

\section{Centrally symmetric polytopes with many edges}
In this section we provide a construction of a cs 2-neighborly
polytope of dimension $d$ and with about
$3^{d/2}\approx (1.73)^d$ vertices as well a construction of
a cs $d$-polytope with $N$ vertices
(for an arbitrarily large $N$) that has about
$\left(1-3^{-d/2}\right)\binom{N}{2}\approx
\left(1-0.58^d\right)\binom{N}{2}$ edges.
Our construction is a slight modification of the one
from \cite{BLN-many}; however our new  trick allows us
to halve the dimension of the polytope from \cite{BLN-many}
while keeping the number of vertices almost the same
as before.

For an integer $m\geq 1$, consider the curve
\begin{equation}  \label{Phi}
\Phi_m: \S \longrightarrow \R^{2(m+1)},\quad \mbox{where} \quad
\Phi_m(t):= \bigl( \cos t, \sin t, \cos 3t, \sin 3t, \ldots,
\cos(3^{m}t), \sin\left(3^{m}t\right)\bigr).
\end{equation}
Note that $\Phi_1=U_2$, see eq.~(\ref{U_k}).
The key to our construction is the following observation.

\begin{lemma}  \label{key}
For an integer $m\geq 1$ and a finite set
$C\subset \S$, define
$$P(C,m)=\conv\left(\Phi_m(t): \, t\in C\right).$$
Then $P(C,m)$ is a polytope of dimension at most
$2(m+1)$ that has $|C|$ vertices.
Moreover, if the elements of $C$ satisfy
\begin{equation}  \label{injective}
3^i t_1 \not\equiv 3^i t_2 \mod 2\pi \quad
\mbox{for all } t_1,t_2\in C \mbox{ such that } t_1\neq t_2,
\, \mbox{ and all }i=1,2,\ldots, m-1,
\end{equation}
then for every pair of distinct points $t_1,t_2 \in C$
that lie on an open arc of length $\pi(1-\frac{1}{3^m})$,
the interval
$[\Phi_m(t_1), \Phi_m(t_2)]$ is an edge of $P(C,m).$
\end{lemma}
\proof To show that $P(C,m)$ has $|C|$ vertices, we
consider the projection
$\R^{2(m+1)} \longrightarrow \R^4$ that forgets all
but the first four coordinates. Since $\Phi_1=U_2$,
the image of $P(C,m)$ is the polytope
\[P(C,1)=\conv\left( U_2(t): \, t \in C \right). \]
By Theorem \ref{pi/2}, the polytope $P(C, 1)$ has $|C|$
distinct vertices: $U_2(t)$ for $t\in C$. Furthermore,
the inverse image of each vertex $U_2(t)$ of $C(m,1)$
in $P(C,m)$ consists of a single vertex $\Phi_m(t)$ of
$P(C,m)$. Therefore, $\Phi_m(t)$ for $t \in C$ are all
the vertices of $P_m$ without duplicates.

To prove the statement about edges, we proceed by induction
on $m$. As $\Phi_1=U_2$, the $m=1$ case follows from \cite{Sm85}
(see Theorem \ref{pi/2} above and the sentence following it).

Suppose now that $m \geq 2$. Let $t_1, t_2$ be two distinct
elements of $C$ that lie on an open arc of length
$\pi(1-\frac{1}{3^m})$. There are two cases to consider.

\medskip \noindent{\em Case I:}
$\, t_1, t_2$ lie on an open arc of length $2\pi/3$.
In this case, the above projection  of $\R^{2(m+1)}$ onto $\R^4$
maps $P(C,m)$ onto $P(C,1)$, and according to the base of induction,
$[\Phi_1(t_1), \Phi_2(t_2)]$ is an edge of $P(C,1)$.
Since the inverse image of a vertex $\Phi_1(t)$ of $P(C, 1)$
in $P(C,m)$ consists of a single vertex $\Phi_m(t)$ of
$P(C,m)$, we conclude that
$\left[ \Phi_m\left(t_1\right), \ \Phi_m\left(t_2\right) \right]$
is an edge of $P(C,m)$.

\medskip \noindent{\em Case II:}
$\, t_1, t_2$ lie on an open arc of length $\pi(1-\frac{1}{3^m})$,
but not on an arc of length $2\pi/3$. (Observe that since 
$3t_1\not\equiv 3t_2 \mod 2\pi$, the points $t_1$ and $t_2$
may not form an arc of length exactly $2\pi/3$.) 
Then $3t_1$ and $3t_2$ do not coincide and lie on an open arc of length
$\pi(1-\frac{1}{3^{m-1}})$.
Consider the projection of $\R^{2(m+1)}$ onto $\R^{2m}$
that forgets the first two coordinates.
The image of  $P(C,m)$ under this projection is
\[ P(3C, m-1), \quad \mbox{where} \quad
3C := \{3t\mod 2\pi : \,\, t\in C\}\subset \S, \]
and since the pair $(3C, m-1)$ satisfies eq.~(\ref{injective}),
by the induction hypothesis, the interval
$$\left[ \Phi_{m-1}\left(3 t_1\right),\
\Phi_{m-1}\left(3 t_2\right) \right]$$
is an edge of $P(3C, m-1)$. By eq.~(\ref{injective}),
the inverse image of a vertex
$\Phi_{m-1}(3t)$ of $P(3C, m-1)$ in $P(C,m)$ consists
of a single vertex $\Phi_m(t)$ of $P(C,m)$,
and hence we infer that
$\left[ \Phi_m\left(t_1\right), \
\Phi_m\left(t_2\right) \right]$
is an edge of $P(C,m)$.
\endproof

We are now in a position to state and prove the main result
of this section. We follow the notation of Lemma \ref{key}.

\begin{theorem} \label{2-neighb}
Fix integers $m\geq 2$ and $s\geq 2$.
Let $A_m \subset \S$ be the set of $2(3^{m}-1)$
equally spaced points:
\[
A_m=\left\{ \frac{\pi (j-1)}{3^{m}-1} :
\quad j=1, \ldots, 2(3^{m}-1)\right\},
\]
and let $A_{m,s}\subset \S$ be the set of
$2(3^{m}-1)$ clusters of $s$ points each,
chosen in such a way that for all $j=1, \ldots, 2(3^{m}-1)$,
the $j$-th cluster lies on an arc of length $10^{-m}$ that contains
the point $\frac{\pi (j-1)}{3^{m}-1}$,
and the entire set $A_{m,s}$ is centrally symmetric. Then
\begin{enumerate}
\item The polytope $P(A_m, m)$ is a centrally symmetric
2-neighborly polytope of dimension $2(m+1)$ that has
$2(3^{m}-1)$ vertices.
\item The polytope $P(A_{m,s},m)$
is a centrally symmetric $2(m+1)$-dimensional polytope
that has $N:=2s(3^{m}-1)$ vertices and at least
$N(N-s-1)/2>\left(1-3^{-m}\right)\binom{N}{2}$ edges.
\end{enumerate}
 \end{theorem}

\proof To see that $P(A_m, m)$ is centrally symmetric,
note that the transformation
$$t\mapsto t+\pi \, \mod 2\pi$$
maps $A_m$ onto itself and also that $\Phi_m(t+\pi)=-\Phi_m(t)$.
The same argument applies to $P(A_{m,s},m)$.

We now show that the dimension of $P(A_m, m)$ is $2(m+1)$.
If not, then the points $\Phi_m(t): t\in A_m$ are all in an
affine  hyperplane
in $\R^{2(m+1)}$, and hence the $2(3^m-1)$ elements of $A_m$
are roots of a trigonometric polynomial of the form
\[
f(t)=c+\sum_{j=0}^m a_j\cos (3^j t)+ \sum_{j=0}^m b_j\sin (3^j t).
\]
Moreover, $a_m$ and $b_m$ cannot both be zero as by our assumption
$f(t)$ has at least $2(3^m-1)$ roots, and so the degree of
$f(t)$ is at least $3^m-1>3^{m-1}$. Thus the complex polynomial
$\P(f)$ defined by eq.~(\ref{P(f)}) is
of the form
\[
\P(f)(z)= d_m z^{2\cdot 3^m} + d_{m-1} z^{3^m+3^{m-1}} +
d_{m-2}z^{3^m+3^{m-2}}+ \cdots +cz^{3^m} + \cdots + \overline{d_m},
\quad \mbox{where } d_m\neq 0.
\]

Note that since $m>1$, $3^m+3^{m-1} < 2\cdot 3^m - 2$.
In particular, the coefficients of $z^{2\cdot 3^m - 1}$
and $z^{2\cdot 3^m - 2}$ are both equal to 0.
Therefore, the sum of all the roots (counted
with multiplicities) of $\P(f)$
as well as the sum of their squares is 0. As
$\deg\P(f)=2\cdot 3^m$, the (multi)set of
roots of $\P(f)$ consists of
$\{e^{it}: \, t\in A_m\}$ together with two additional roots,
denote them by $\zeta_1$ and $\zeta_2$. The complex numbers
$e^{it}: t\in A_m$ form a geometric progression,
and it is straightforward to check that
\[
\sum_{t\in A_m} e^{it}=0 \quad \mbox{and} \quad
\sum_{t\in A_m} e^{2it}=0.
\]
Hence for the sum of all the roots of $\P(f)$
and for the sum of their squares to be zero, we must
have
\[\zeta_1+\zeta_2=0 \quad \mbox{and} \quad
\zeta_1^2+\zeta_2^2=0.
\]
Thus $\zeta_1=\zeta_2=0$, and so the constant term of
$\P(f)$ is zero. This however contradicts the fact that
the constant term of $\P(f)$ equals $\overline{d_m}$,
where $d_m\neq 0$. Therefore, the polytope $P(A_m,m)$ is
full-dimensional.

Finally, to see that $P(A_m,m)$ is 2-neighborly,
observe that it follows from the definition of $A_m$ that if
$t_1, t_2\in A_m$ are not antipodes, then they lie on a
closed arc of length $\pi(1-\frac{1}{3^m-1})$, and
hence also on an an open arc of length $\pi(1-\frac{1}{3^m})$.
In addition, since $3^m-1$ is relatively prime to 3,
we obtain that for every two distinct elements $t_1$, $t_2$ of $A_m$,
$3^i t_1 \not\equiv 3^i t_2 \mod 2\pi$ (for $i=1,\ldots,m-1$).
Part (1) of the theorem is then immediate from Lemma~\ref{key}.

To compute the dimension of $P(A_{m,s},m)$, note that if
it is smaller than $2(m+1)$, then $P(A_{m,s},m)$ is a subset of an
affine hyperplane in $\R^{2(m+1)}$. As all vertices of this polytope
lie on the curve $\Phi_m$, such a hyperplane corresponds to a
trigonometric polynomial of degree $3^m$ that has
at least $N=2s(3^m-1)\geq 4(3^m-1) > 2\cdot 3^m$ roots.
This is however impossible, as no nonzero trigonometric polynomial
of degree $D$ has more than $2D$ roots.

To finish the proof of Part (2),
note that since each cluster of $A_{m,s}$ lies
on an open arc of length 
$$10^{-m}<\frac{\pi}{2}\left(\frac{1}{3^m-1}-\frac{1}{3^m}\right)$$ 
that contains the corresponding element
of $A_m$, and since multiplication by $3^i$ modulo $2\pi$
maps $A_m$ bijectively onto itself, it follows that
\begin{itemize}
\item $3^i t_1 \not\equiv 3^i t_2 \mod 2\pi$
(for $i=1,\ldots,m-1$) holds for all distinct
$t_1, t_2\in A_{m,s}$. 
(Indeed, for $t_1$, $t_2$ from the same cluster, 
the points $3^it_1$ and $3^it_2$ of $\S$ 
do not coincide as $3^m/10^m < 2\pi$, 
and for $t_1$, $t_2$ from different clusters, 
$3^it_1$ and $3^it_2$ do not coincide as
the distance between them
along $\S$ is at least 
$\frac{\pi}{3^m-1}-\frac{2\cdot 3^m}{10^m}>0$.)

\item Every two points $t_1, t_2\in A_{m,s}$
lie on an open arc of length $\pi(1-\frac{1}{3^m})$
as long as they do not belong to a pair of
opposite clusters.
\end{itemize}
Thus  Lemma \ref{key} applies and shows that
 the interval $\left[\Phi_m(t_1), \,  \Phi_m(t_2)\right]$
is an edge of $P(A_{m,s},m)$ for all $t_1,t_2\in A_{m,s}$
that are not from opposite clusters. In other words,
each vertex of $P(A_{m,s},m)$ is incident with at
least $N-s-1$ edges. This yields the promised bound on the number
of edges of $P(A_{m,s},m)$ and completes the proof of Part (2).
\endproof

\section{Applications to strict antipodality problems}
In this section we observe that an appropriately chosen half
of the vertex set of any cs $2k$-neighborly
 $d$-dimensional polytope has a large number of
pairwise strictly antipodal $(k-1)$-simplices.
The results of the previous section then imply new lower bounds
on questions related to strict antipodality. Specifically,
in the following theorem we improve both Talata's and
Makai--Martini's bounds.

\begin{theorem} \label{A'_d} \quad
\begin{enumerate}
\item
For every $m\geq 1$, there exists a set
$X_m \subset \R^{2(m+1)}$ of size $3^m-1$ that affinely spans
$\R^{2(m+1)}$ and such that each pair
of points of $X_m$ is strictly antipodal. Thus,
$A'(d)\geq 3^{\lfloor d/2-1\rfloor}-1$ for all $d\geq 4$.
\item
For all positive integers $m$ and $s$, there exists a set
$Y_{m,s}\subset\R^{2(m+1)}$ of size $n:=s(3^m-1)$ that has
at least
\[
\left( 1-\frac{1}{3^m-1}\right)\cdot \frac{n^2}{2}
\]
pairs of antipodal points. Thus,
$A'_d(n)\geq
\left(1-\frac{1}{3^{\lfloor d/2-1\rfloor}-1}\right)\cdot
\frac{n^2}{2}-O(n)$
for all $d\geq 4$ and $n$.
\end{enumerate}
\end{theorem}

One can generalize the notion of strictly antipodal points
in the following way: for a set $X\subset \R^d$ that affinely
spans $\R^d$, we say that two simplices,
$\sigma$ and $\sigma'$, spanned by the points of $X$
 are strictly antipodal if there exist two distinct parallel
hyperplanes $H$ and $H'$ such that $X$ lies in the slab defined
by $H$ and $H'$,
$H\cap \conv(X)=\sigma$, and $H'\cap \conv(X)=\sigma'$.
Makai and Martini \cite{MakMart} asked about the maximum number
of pairwise strictly antipodal $(k-1)$-simplices in $\R^d$.
The following result gives a lower bound to their question.

\begin{theorem} \label{simplices}
There exists a set of $\lfloor(d/2)\cdot 2^{cd/k}\rfloor$
points in $\R^d$ with the property that every two disjoint
$k$-subsets of $X$ form the vertex sets of strictly antipodal
$(k-1)$-simplices. In particular, there exists a set of
$\lfloor \frac{d}{2k}\cdot 2^{cd/k}\rfloor$
pairwise strictly antipodal $(k-1)$-simplices in $\R^d$.
Here $c>0$ is an absolute constant.
\end{theorem}

The key to our proofs are results of Section 3 and
paper \cite{LinNov} along with the following observation.

\begin{lemma} \label{half}
Let $P\subset \R^d$ be a full-dimensional cs polytope on the
vertex set $V=X\sqcup (-X)$. If $U_1$, $U_2$ are
subsets of $X$ such that $U_1\cup (-U_2)$ is the vertex set
of a $(|U_1|+|U_2|-1)$-face of $P$, then $\sigma_1:=\conv(U_1)$
and $\sigma_2:=\conv(U_2)$ are strictly antipodal simplices
spanned by points of $X$. In particular, if $P$ is $2$-neighborly,
then every pair of vertices of $X$ is strictly antipodal, and,
more generally, if $P$ is $2k$-neighborly, then every two disjoint
$k$-subsets of $X$ form a pair of strictly antipodal $(k-1)$-simplices.
\end{lemma}
\proof
Since $\tau_1:=\conv(U_1\cup(-U_2))$ is a face of $P$,
there exists a supporting hyperplane
$H_1$ of $P$ that defines $\tau_1$:
specifically, $P$ is contained in one of the closed half-spaces
bounded by $H_1$ and $P\cap H_1=\tau_1$. As $P$ is centrally symmetric,
the hyperplane $H_2:=-H_1=\{x\in \R^d \ : -x\in H_1\}$
is a supporting hyperplane of $P$ that defines
the opposite face, $\tau_2:=\conv((-U_1)\cup U_2)$.
Thus $P$, and hence also $X$, is contained in the slab
between $H_1$ and $H_2$. Moreover, since $U_1, U_2$ are
subsets of $X$, it follows that $-U_1$ and $-U_2$
are contained in $-X$, and hence disjoint from $X$.
Therefore,
\[H_i\cap \conv(X) =H_i\cap P\cap \conv(X)=
\tau_i\cap \conv(X)=\conv(U_i)=\sigma_i \quad \mbox{for $i=1,2$}.
\]
The result follows.
\endproof

\smallskip\noindent{\em Proof of Theorem \ref{A'_d}: \, }
Consider the sets $A_m$ and $A_{m,s}$ of Theorem~\ref{2-neighb}.
Define $$A_m^+=\{t\in A_m : \, 0\leq t <\pi\},$$
and define $A_{m,s}^+$ by taking the union of
those clusters of $A_{m,s}$ that lie on small arcs
around the points of $A_m^+$. In particular,
$|A_m^+|=3^m-1$ and $|A_{m,s}^+|=s(3^m-1)$.
Let
$$X_m:=\left\{\Phi_m(t) :
\, t\in A_m^+\right\}\subset\R^{2(m+1)} \quad \mbox{and} \quad
Y_{m,s}:=\left\{\Phi_m(t) :
\, t\in A_{m,s}^+\right\}\subset\R^{2(m+1)}.$$
Theorem \ref{2-neighb} and Lemma \ref{half} imply that
each pair of points of $X_m$ is strictly antipodal,
and each pair of points of $Y_{m,s}$ that are not from the
same cluster is strictly antipodal. The claim follows.
\endproof

\smallskip\noindent{\em Proof of Theorem \ref{simplices}: \, }
It was proved in \cite{LinNov, RV} (by using a probabilistic
construction) that if $k$, $d$, and $N$ satisfy
\[
k\leq\frac{cd}{1+\log \frac{N}{d}},
\]
where $c>0$ is some absolute constant, then there exists a
$d$-dimensional cs polytope on $2N$ vertices that is
$2k$-neighborly.  Solving this inequality for $N$,
implies existence of a $d$-dimensional cs polytope on
$\lfloor d\cdot 2^{cd/k}\rfloor $ vertices that is $2k$-neighborly.
This together with Lemma \ref{half} yields the result.
\endproof

\section{Constructing $k$-neighborly cs polytopes}
The goal of this section is to present a deterministic
construction of a cs $k$-neighborly
$d$-polytope with at least $2^{c_k d}$ for $c_k=3/20k^2 2^k$ vertices.
This requires the following facts and definitions.

A family $\F$ of subsets of $[m]:=\{1,2,\ldots,m\}$
is called {\em $k$-independent} if for every $k$ distinct
subsets $I_1,\ldots, I_k$ of $\F$ all $2^k$ intersections
$$\bigcap_{j=1}^k J_j, \quad \mbox{where } J_j=I_j \mbox{ or }
J_j= I_j^c:=[m]\setminus I_j,  \,\mbox{ are non-empty}.$$
The crucial component of our construction is a
deterministic construction of $k$-independent families of size
larger than $2^{m/5(k-1)2^k}$ given in \cite{FLL}.

For a subset $I$ of $[m]$ and a given number $a\in\{0,1\}$,
we (recursively) define a sequence
$x(I, a)=(x_0, x_1,\ldots, x_m)$
of {\bf zeros and ones} according to the following rule:
\begin{equation}  \label{x(I,alpha)}
x_0=x_0(I, a):=a \quad \mbox{and} \quad
x_{n}=x_n(I, a) \equiv
\left\{ \begin{array}{ll}
    \sum_{j=0}^{n-1} x_j & \mbox{ if $n\notin I$}\\
  1+ \sum_{j=0}^{n-1} x_j & \mbox{ if $n\in I$}
 \end{array}
\right.
\mod 2 \quad \mbox{for $n\geq 1$}.
\end{equation}
We also set
\begin{equation} \label{t(I,alpha)}
t(I,a):=\pi \sum_{j=0}^m \frac{x_j}{3^j} \in \S.
\end{equation}

A few observations are in order.
First, it follows from (\ref{x(I,alpha)}) that
$x(I,a)\neq x(J,a)$ if $I\neq J$, and that
 $x(I, a)$ and $x(I^c, 1-a)$ agree in
all but the 0-th component, where they disagree. Hence
$$t(I, a)=t(I^c, 1-a)+\pi \mod 2\pi.$$
Second, since $\sum_{j=1}^\infty \frac{1}{3^j}=\frac{1}{2}$
and since all components of $x(I, a)$ are zeros and ones,
we infer from eq.~(\ref{t(I,alpha)})
that for all $1\leq n \leq m$ and all
$0\leq \epsilon \leq 1/3^{m+1}$, the point
$3^n\cdot\left(t(I,a)+\pi\epsilon\right)$ of $\S$ either
lies on the arc $[0,\pi/2)$ or on the arc $[\pi,3\pi/2)$
depending on the parity of
$$\sum_{j=0}^n 3^{n-j} x_j(I,a)
\equiv \sum_{j=0}^n x_j(I,a)
\mod 2 .$$
As, by (\ref{x(I,alpha)}), $\sum_{j=0}^n x_j(I,a)$ is even if
$n\notin I$ and is odd if $n\in I$, we obtain that
\begin{equation} \label{pi,3pi/2}
3^n\cdot\left(t(I,a)+\pi\epsilon\right)
\in [\pi, 3\pi/2) \mod 2\pi \quad \
\mbox{for all $n\in I$ and $a\in\{0,1\}$}.
\end{equation}

The relevance of $k$-independent sets to cs $k$-neighborly polytopes
is explained by the following lemma along with Theorem \ref{pi/2}.

\begin{lemma}   \label{small-arc}
Let $\F$ be a $k$-independent family of subsets of $[m]$, let
$\epsilon_I\in[0, 1/3^{m+1}]$ for $I\in \F$, and let
$$V^\epsilon(\F)= \bigcup_{I\in\F} \left\{
t\left(I,0\right)+\pi\epsilon_I, \,
t\left(I^c,1\right)+\pi\epsilon_I\right\}
\subset \S.$$
Then for every $k$ distinct points
$t_1,\ldots, t_k$ of $V^\epsilon(\F)$ no two of which
are antipodes, there exists an integer $n\in[m]$ such that
the subset $\{3^n t_1, \ldots, 3^nt_k\}$ of $\S$
is entirely contained in $[\pi,3\pi/2)$.
\end{lemma}
\proof As $t_1,\ldots, t_k$ are elements of $V^\epsilon(\F)$,
by relabeling them if necessary,
we can assume that
$$
t_j=\left\{\begin{array}{ll} t(I_j, 0)+\pi\epsilon_{I_j} &
\mbox{ if $1\leq j \leq q$}\\
t(I_j^c, 1)+\pi\epsilon_{I_j} & \mbox{ if $q<j\leq k$}
\end{array}
\right.
$$
for some $0\leq q\leq k$ and $I_1,\ldots, I_k\in\F$.
Moreover, the sets $I_1,\ldots, I_k$ are distinct,
since $t_1,\ldots, t_k$ are distinct and no two of them
are antipodes.
As $\F$ is a $k$-independent family, the intersection
$(\cap_{j=1}^q I_j)\cap (\cap_{j=q+1}^k I_j^c)$ is non-empty.
The result follows, since by eq.~(\ref{pi,3pi/2}),
for any element $n$ of this intersection,
$\{3^nt_1,\ldots, 3^nt_k\}\subset[\pi,3\pi/2)$.
\endproof

For $I\in \F$, define $\epsilon_I=\epsilon_{I^c}
:= \sum_{i\in I} 10^{-i-m}.$
Then 
\begin{equation} \label{injectivite'}
3^n t_1\not\equiv 3^n t_2 \mod 2\pi \quad \mbox{for all }
t_1,t_2\in V^\epsilon(\F) \mbox{ such that } t_1\neq t_2,
\, \mbox{ and all } 1\leq n \leq m.
\end{equation}  
Indeed, if $t_1$ and $t_2$ are antipodes, then so are $3^n t_1$
and $3^n t_2$, and (\ref{injectivite'}) follows. If $t_1$
and $t_2$ are not antipodes, then
there exist two distinct and not complementary subsets 
$I, J$ of $[m]$ such that $t_1=t(I, a)+\pi\epsilon_I$
and $t_2=t(J,b)+\pi\epsilon_J$ for some $a,b\in\{0,1\}$.
Hence, by  definition of $\epsilon_I$ and $\epsilon_J$,
\[ \pi/10^{2m}  <
  3^n\cdot \pi|\epsilon_I-\epsilon_J| <
  \pi (3/10)^{m},
\]
while by definition of $t(I,a)$ and $t(J,b)$,
the distance between the points 
$3^n \cdot t(I,a)$ and $3^n \cdot t(J,b)$ of $\S$ 
along $\S$ is either $0$ or at least $\pi/3^m$.
In either case, it follows that the distance
between $3^n\left(t(I, a)+\pi\epsilon_I\right)$
and $3^n\left(t(J, b)+\pi\epsilon_J\right)$
is positive, yielding eq.~(\ref{injectivite'}).

We are now in a position to present our construction of
$k$-neighborly cs polytopes.  The construction is similar to
that in Theorem \ref{2-neighb}, except that it is based on
the set $V^\epsilon(\F)\subset \S$, where $\F$ is a
$k$-independent family of subsets of $[m]$, instead of
$A_m\subset \S$, and on a modification of $\Phi_m$ to
a curve that involves $U_k$ instead of $U_2$.

Let $U_k: \S\longrightarrow\R^{2k}$ be the curve defined by eq.~(\ref{U_k}).
In analogy with the curve $\Phi_m$ (see eq.~(\ref{Phi})),
for integers $m\geq 0$ and $k\geq 3$, define the curve
\begin{equation}   \label{Psi}
\Psi_{k,m}: \S\longrightarrow\R^{2k(m+1)} \quad \mbox{by} \quad
\Psi_{k,m}(t):=\left(U_k(t), U_k(3t), U_k(3^2t),
\ldots, U_k(3^m t)\right).
\end{equation}
Thus, $\Psi_{k,0}=U_k$ and $\Psi_{k,m}(t+\pi)=-\Psi_{k,m}(t)$.

The following theorem is the main result of this section.
We use the same notation as in Lemma~\ref{small-arc}.
Also, mimicking the notation of Lemma \ref{key},
for a subset $C$ of $\S$, we denote by $P_k(C,m)$ the polytope
$\conv(\Psi_{k,m}(t): \,\, t\in C)$.

\begin{theorem} \label{k-neighb}
Let $m\geq 1$ and $k\geq 3$ be fixed integers,
let $\F$ be a $k$-independent family of subsets of $[m]$,
and let $\epsilon_I=\sum_{i\in I} 10^{-i-m}$ for $I\in \F$.
Then the polytope
\[P_{k}(V^\epsilon(\F), m):=\conv\left(\Psi_{k,m}(t) :
\quad t\in V^\epsilon(\F)\right) \]
is a cs $k$-neighborly polytope of dimension
at most $2k(m+1)-2m\lfloor(k+1)/3\rfloor$ that has
$2|\F|$ vertices.
\end{theorem}

\begin{remark} \label{numerics}
For a fixed $k$ and an arbitrarily large $m$,
a deterministic algorithm from  \cite{FLL} produces
a $k$-independent family $\F$ of subsets of $[m]$
such that $|\F|>2^{m/5(k-1)2^k}$. Combining this
with Theorem \ref{k-neighb} results in a
 cs neighborly polytope of dimension $d\approx \frac{4}{3}km$
and more than $2^{{3d}/{20k^2 2^k}}$ vertices. Of a special interest
is the case of $k=3$: the algorithm from \cite{FLL} provides
a 3-independent family of size $\approx 2^{0.092m}$, which together
with Theorem \ref{k-neighb} yields a deterministic
construction of a cs 3-neighborly polytope of dimension $\leq d$
and with about $2^{0.023d}$ vertices.
\end{remark}

\smallskip\noindent{\it Proof of Theorem \ref{k-neighb}:\,}
As in the proof of Theorem \ref{2-neighb},
the polytope $P_{k}(V^\epsilon(\F), m)$ is centrally symmetric
  since $V^\epsilon(\F)$ is a cs subset
of $\S$ and since $\Psi_{k,m}(t+\pi)=-\Psi_{k,m}(t)$.

Also as in the proof of Theorem \ref{2-neighb}, the fact
that $P_{k}(V^\epsilon(\F), m)$ has $2|\F|$ vertices follows by
considering the projection $\R^{2k(m+1)}\longrightarrow \R^{2k}$ that
forgets all but the first $2k$ coordinates. Indeed, the
image of $P_{k}(V^\epsilon(\F), m)$ under this projection is
the polytope
$$P_{k}(V^\epsilon(\F), 0)=
\conv\left(U_k(t): \, t\in V^\epsilon(\F)\right),$$
and this latter polytope has $2|\F|$ vertices
(by Theorem \ref{U_k}).

To prove $k$-neighborliness of $P_{k}(V^{\epsilon}(\F), m)$,
let $t_1,\ldots,t_k\in V^{\epsilon}(\F)$ be $k$
distinct points no two of which are antipodes.
By Lemma \ref{small-arc}, there exists an integer
$1\leq n\leq m$ such that the points $3^nt_1,\ldots, 3^nt_k$
of $\S$ are all contained in the arc $[\pi, 3\pi/2)$.
Consider the projection $\R^{2k(m+1)}\longrightarrow \R^{2k(m+1-n)}$
that forgets the first $2kn$ coordinates followed by the
projection $\R^{2k(m+1-n)} \longrightarrow \R^{2k}$ that forgets all but
the first $2k$ coordinates. The image of
$P_{k}(V^\epsilon(\F), m)$ under this composite projection
is
$$P_{k}(3^n V^{\epsilon}(\F), 0)
=\conv (U_k(3^n t) : \quad t\in  V^{\epsilon}(\F)),
$$
and, since $\{3^nt_1,\ldots, 3^nt_k\}\subset[\pi,3\pi/2)$,
Theorem \ref{U_k} implies that
the set $\{U_k(3^n t_i) : i=1,\ldots,k\}$
is the vertex set of a $(k-1)$-face of this latter
polytope. As, by eq.~(\ref{injectivite'}), the inverse
image of a vertex $U_k(3^nt)$ of
$P_{k}(3^n V^{\epsilon}(\F), 0)$
in $P_{k}(V^{\epsilon}(\F), m)$ consists of a single
vertex $\Psi_{k,m}(t)$ of $P_{k}(V^{\epsilon}(\F), m)$,
we obtain that $\{\Psi_{k,m}(t_i) : i=1,\ldots,k\}$
is the vertex set of a $(k-1)$-face of
$P_{k}(V^{\epsilon}(\F), m)$. This completes
the proof of $k$-neighborliness of
$P_{k}(V^{\epsilon}(\F), m)$.

To bound the dimension of $P_{k}(V^{\epsilon}(\F), m)$,
observe that the third  coordinate of $U_k(t)$ coincides
with the first coordinate of $U_k(3t)$ while the fourth
coordinate of $U_k(t)$ coincides with the second coordinate
of $U_k(3t)$, etc. Thus $P_{k}(V^{\epsilon}(\F), m)$ is in a
subspace of $\R^{2k(m+1)}$,
and to bound the dimension of this subspace
we must account for all repeated coordinates. This
can be done exactly as in \cite[Lemma~2.3]{BLN-many}.
We leave details to our readers.
\endproof

Fix $s\geq 2$, and let $V^{\epsilon}(\F,s)$ be a centrally
symmetric subset of $\S$ obtained by replacing each point
$t\in V^\epsilon(\F)$ (in Theorem \ref{k-neighb}) with a
cluster of $s$ points that all lie on a sufficiently
small open arc containing $t$. Then the proof
of Theorem \ref{k-neighb} implies that the polytope
$P_k(V^{\epsilon}(\F,s),m)$ is a cs polytope
with $N:=2s|\F|$ vertices, of
dimension at most $2k(m+1)-2m\lfloor(k+1)/3\rfloor$,
and such that every $k$
vertices of this polytope no two of which are from opposite
clusters form the vertex set of a $(k-1)$-face.
Choose a $k$-element set from the union of these
$2|\F|$ clusters (of $s$ points each) at random from
the uniform distribution.
Then the probability that this set has no two points from
opposite clusters is at least
\[ \prod_{i=0}^{k-1} \frac{(2|\F|-i)s-i}{2|\F|s-i}
\geq \prod_{i=0}^{k-1}
\left(1-\frac{i}{|\F|}\right)
\geq 1-\frac{k^2}{|\F|}.
\]
Thus, the resulting polytope has at least
 $$\left(1-\frac{k^2}{|\F|}\right)\binom{N}{k}$$
$(k-1)$-faces.
Combining this estimate with Remark \ref{numerics},
we obtain

\begin{corollary}  \label{many-k}
For a fixed $k$ and arbitrarily large $N$ and $d$, there exists
a cs $d$-dimensional polytope with $N$ vertices and at least
\[
\left(1- k^2
\left(2^{-{3}/{20k^2 2^k}}\right)^d \right)\binom{N}{k}
\]
$(k-1)$-faces.
\end{corollary}
This corollary  improves  \cite[Cor.~1,4]{BLN-many}
asserting existence of cs $d$-polytopes with $N$ vertices
and at least $\left(1-(\delta_k)^d\right)\binom{N}{k}$
faces of dimension $k-1$, where
$\delta_k\approx (1-5^{-k+1})^{5/(24k+4)}$.

\section*{Acknowledgements}
We are grateful to Zolt\'an F\"uredi and Noga Alon for 
pointing us to efficient constructions of 
$k$-independent families and to 
Endre Makai  for references to recent results on antipodal pairs.

\small{
}


\begin{thebibliography}{99}

\bibitem{Alon} N.~Alon, 
``Explicit construction of exponential sized families of 
$k$-independent sets'',
 Discrete Math.~{\bf 58} (1986), no.~2, 191--193.

\bibitem{AMS} N.~Alon, D.~Moshkovitz, and S.~Safra, 
``Algorithmic construction of sets for $k$-restrictions'', 
ACM Trans.~Algorithms {\bf 2} (2006), no.~2, 153--177. 


\bibitem{BN} A.~Barvinok and I.~Novik,
``A centrally symmetric version of the cyclic polytope",
Discrete Comput.~Geometry {\bf 39} (2008), 76--99.

\bibitem{BLN-many}
A.~Barvinok, S.~J.~Lee, and I.~Novik,
``Centrally symmetric polytopes with many faces'',
Israel J.~Math., to appear, arXiv:1106.0449 [math.MG].

\bibitem{BLN-neighb} A.~Barvinok, S.~J.~Lee, and I.~Novik,
``Neighborliness of the symmetric moment curve'',
Mathematika, to appear, doi:10.1112/S0025579312000010.

\bibitem{Bor} K.~B\"or\"oczky, Jr.,
{\em Finite packing and covering},
Cambridge Tracts in Mathematics,
154, Cambridge University Press, Cambridge, 2004.

\bibitem{DanzGr} L.~Danzer and B.~Gr\"unbaum,
``\"Uber zwei Probleme bez\"uglich konvexer
K\"orper von P.~Erd\"os und von V.~L.~Klee''
(German), Math.~Z. {\bf 79} (1962) 95--99.

\bibitem{ErdFur}
P.~Erd\H{o}s and Z.~F\"uredi, ``The greatest angle
among $n$ points in the $d$-dimensional Euclidean space'',
Combinatorial mathematics (Marseille-Luminy, 1981), 275--283,
North-Holland Math. Stud., 75, North-Holland, Amsterdam, 1983.

\bibitem{FLL} G.~Freiman, E.~Lipkin, and L.~Levitin,
``A polynomial algorithm for constructing families of
$k$-independent sets'', Discrete Math.~{\bf 70} (1988),
137-147.

\bibitem{Ha} V.~Harangi, 
``Acute sets in Euclidean spaces'', 
SIAM J.~Discrete Math.~{\bf 25} (2011), no. 3, 1212--1229.

\bibitem{Indyk}  P.~Indyk, 
``Uncertainty principles, extractors, and explicit embeddings of $l_2$ 
into $l_1$'', 
STOC'07-- Proceedings of the 39th Annual ACM Symposium on 
Theory of Computing, 615--620, ACM, New York, 2007.

\bibitem{LinNov}  N.~Linial and I.~Novik,
``How neighborly can a centrally symmetric polytope be?'',
Discrete Comput.~Geometry {\bf 36} (2006), 273--281.


\bibitem{MakMart}  E.~Makai, Jr.~and H.~Martini,
``On the number of antipodal or strictly antipodal
pairs of points in finite subsets of ${\bf R}^d$'',
in: P.~Gritzmann, B.~Sturmfels (Eds.),
Applied geometry and discrete mathematics, 457--470,
DIMACS Ser. Discrete Math. Theoret. Comput. Sci., 4,
Amer. Math. Soc., Providence, RI, 1991.

\bibitem{MartSolt} H.~Martini and V.~Soltan,
``Antipodality properties of finite sets in Euclidean space'',
Discrete Math.~{\bf 290} (2005), no. 2-3, 221--228.

\bibitem{McM-UBC}
P.~McMullen, ``The maximum numbers of faces of a convex polytope'',          
Mathematika {\bf 17} (1970), 179--184.   

\bibitem{RV} M.~Rudelson and R.~Vershynin,
``Geometric approach to error correcting codes and
reconstruction of signals",
 Int.~Math.~Res.~Not., {\bf 64} (2005), 4019--4041.

\bibitem{Sm85} Z.~Smilansky,
``Convex hulls of generalized moment curves",
Israel J.~Math.~{\bf 52} (1985), 115--128.
\end{thebibliography}
\end{document}